\documentclass[cupthm]{CUP-JNL-JMJ}

\usepackage{latexsym}
\usepackage{graphicx}
\usepackage{multicol,multirow}
\usepackage{amsmath,amssymb,amsfonts}
\usepackage{mathrsfs}
\usepackage{amsthm}
\usepackage{rotating}
\usepackage{appendix}
\usepackage[numbers]{natbib}
\usepackage{float} 

\numberwithin{equation}{section}

\theoremstyle{cupplain}
\newtheorem{theorem}{Theorem}[section]
\newtheorem{lemma}[theorem]{Lemma}

\newtheorem{proposition}[theorem]{Proposition}

\theoremstyle{cupdefinition}
\newtheorem{definition}{Definition}[section]

\theoremstyle{cupremark}

\theoremstyle{cupproof}
\newtheorem{proof}{Proof}

\def\N{{\mathbb N}}

\def\R{{\mathbb R}}
\def\Z{{\mathbb Z}}

\def\cF{{\mathcal F}}

\def\cR{{\mathcal R}}

\def\cT{{\mathcal T}}

\def\Spec{{\rm Spec\,}}
\def\sss{{{\mathbb S}}}
\def\spm{{\sss[\pm 1]}}
\def\Se{\frak{ Sets}}
\def\Ses{{\Se_*}}
\def\Hom {{\rm{Hom}}}

\def\spzb{{\overline{\Spec\Z}}}

\def\gop{{\Gamma^{\rm op}}}
\def\spzb{{\overline{\Spec\Z}}}

\newcommand{\ie}{{\it i.e.\/}\ }

\begin{document}

\begin{Frontmatter}

\title[Riemann-Roch for the ring $\Z$]{\bf Riemann-Roch for the ring $\Z$\thanks{Research supported by The Simons Foudation}}

\author[1,2]{ALAIN CONNES}

\author[3]{CATERINA CONSANI}

\address[1]{\orgname{Coll\`ege de France} 
\orgaddress{\city{Paris}, \state{F-75005 France}}

\email{alain@connes.org}

}

\address[2]{\orgname{IHES}
\orgaddress{\city{Bures-sur-Yvette} \state{91440 France}}

}

\address[3]{\orgdiv{Department of Mathematics} \orgname{The Johns Hopkins University}, 
\orgaddress{\city{Baltimore} \state{21218 USA}}

\email{cconsan1@jhu.edu}

}

\maketitle

\authormark{A. Connes and C. Consani}

\abstract{We show that by working over the absolute base $\sss$ (the categorical version of the sphere spectrum) instead of $\spm$  improves our previous  Riemann-Roch formula for $\spzb$.   The formula equates the (integer-valued) Euler characteristic of an Arakelov divisor with  the sum of the degree of the divisor (using logarithms with base 2)  and the  number $1$, thus confirming the understanding of the ring $\Z$ as a ring of polynomials in one variable over the absolute base $\sss$, namely  $\sss[X], 1+1=X+X^2$. }


\keywords[\textup{2010} Mathematics subject classification]{14C40; 14G40; 14H05; 11R56; 13F35; 18N60; 19D55}


\end{Frontmatter}

\section{Introduction}
\label{intro}
In \cite{RR} we proved a Riemann-Roch formula for $\spzb$ working over  the spherical extension $\spm:=\sss[\mu_{2,+}]$ of the absolute base $\sss$. The proof of that result is based on viewing the ring $\Z$  as a ring of polynomials\footnote{More precisely every integer is uniquely of the form $P(X)$ where $P$ is a polynomial with coefficients in $\{-1,0,1\}$ and $X=3$, the  presentation is given by $1+1=X-1$} with coefficients in $\spm$ and generator $3\in \Z$. In the present paper we show  that by working over the absolute base $\sss$ itself,  one obtains the following Riemann-Roch formula.
\begin{theorem}\label{rrspzbintro} 
Let $D$ be an Arakelov divisor  on $\spzb$. Then\footnote{We use the notation $\deg_2:=\deg/\log 2$ } 
\begin{equation}\label{rrforq}
\dim_{\sss}H^0(D)-\dim_{\sss}H^1(D)=\bigg\lceil \deg_2 D\bigg\rceil'	+1.
\end{equation}
Here $\lceil x \rceil'$ denotes the right continuous function which agrees with the function ceiling$(x)$ for $x>0$ non-integer, and with --ceiling$(-x)$ for $x<0$ non-integer (see Figure \ref{rr1}).
\end{theorem}
The proof of \eqref{rrforq} follows the same lines as the proof of the Riemann-Roch formula in \cite{RR}, and views   $\Z$ as a ring of polynomials\footnote{Every integer is uniquely of the form $P(X)$ where $P$ is a polynomial with coefficients in $\{0,1\}$ and $X=-2$, the presentation is $1+1=X+X^2$} over $\sss$ with generator $-2$. It greatly improves this earlier result as follows: 
\begin{enumerate}
\item The term $\mathbf 1_L$ involving the exceptional set $L$ in the earlier formula is now eliminated.
\item Formula \eqref{rrforq} displays a  perfect analogy with the Riemann-Roch formula holding for curves of genus $0$.
\item The canonical divisor  $K=-2\{2\}$ has integral degree $\deg_2(K)=-2$.
\end{enumerate}
When working over the absolute base $\sss$ one is led to a very natural notion of $\sss$-module associated to an Arakelov divisor as explained in Section \ref{sect1}.

\section{Working over the absolute base $\sss$}\label{sect1}
We let $\gop$ be the opposite of the Segal category (see  \cite{DGM} Chpt. 2 and \cite{CCprel}), it has one  object $k_+$ for each integer $k>0$, the  pointed set $\{*,1,\ldots,k\}$,  and the morphisms are morphisms of pointed sets. Covariant functors $\gop\longrightarrow \Se_*$ and their natural transformations determine the category $\Gamma\Ses$ of $\Gamma$-sets (aka $\sss$-modules). 
When working over the spherical monoidal algebra $\spm$ of the (pointed) multiplicative monoid $\{\pm 1\}$, the natural $\spm$-module associated to a norm on an abelian group $A$ is ($k\in\N$, $\lambda\in\R$)
\begin{equation}\label{normdefn}
\Vert HA\Vert_\lambda(k_+):=\{ a\in A^k\mid \sum \vert a_j\vert \leq \lambda\}.
\end{equation}
The above formula  is applied at the archimedean place, for subgroups $A\subset \R$ and with $\vert \cdot \vert$ denoting the euclidean absolute value. If $\spm$ is replaced by the base $\sss$, there is a more basic definition of an $\sss$-module associated to an arbitrary subset $X\subset A$ containing $0\in A$ 

\begin{lemma}\label{xlem} Let $A$ be an abelian monoid with $0\in A$. Let $X\subset A$ be a subset containing $0$. The following condition defines a subfunctor of the $\sss$-module $HA$ 
	\begin{equation}\label{xdefn}
(HA)_X(k_+):=\{ a\in A^k\mid \sum_Z a_j \in X, \  \forall Z\subset k_+\}\subset X^k.
\end{equation}
\end{lemma}
\proof By construction $(HA)_X(k_+)$ is a subset of 
$HA(k_+)$ containing the base point $a_j=0$, $\forall j$. Let $\phi:k_+\to m_+$ be a map preserving the base point $*$, we shall show that $\phi_*((HA)_X(k_+))\subset (HA)_X(m_+)$. Let $a \in (HA)_X(k_+)$. For any $\ell \in m_+$, $\ell\neq  *$, one has 
$$
\phi_*(a)(\ell)=\sum_{\phi^{-1}(\ell)} a_j=\sum_{Z_\ell}a_j, \qquad Z_\ell:=\phi^{-1}(\ell).
$$
It follows from \eqref{xdefn} that $\phi_*(a)(\ell)\in X$ for all $\ell$ and  that for any pointed subset $Z'\subset m_+$ 
$$
\sum_{\ell \in Z'} \phi_*(a)(\ell)=\sum_{Z}a_j \in X,  \qquad Z=\cup_{\ell \in Z'} Z_\ell.
$$
This proves that $\phi_*((HA)_X(k_+))\subset (HA)_X(m_+)$.\endproof 
Next Proposition shows that for $X= [-\lambda,\lambda]\subset \R$ a symmetric interval, the $\sss$-module $(H\R)_X$ is a module over the $\sss$-algebra $\Vert H\R\Vert_1$.
\begin{proposition}\label{higherspm} Let $\lambda>0$,  $X= [-\lambda,\lambda]\subset \R$ a symmetric interval and $(H\R)_X$ as in \eqref{xdefn}. Then
\begin{equation}\label{xdefnr}
(H\R)_X(k_+)=\{ a\in \R^k\mid \sum_{a_j> 0} a_j \leq \lambda, \   \sum_{a_j<0} (-a_j) \leq \lambda\}
\end{equation}
Moreover, the module action  of the $\sss$-algebra $H\R$ on itself by multiplication  induces an action of the $\sss$-algebra $\Vert H\R\Vert_1$ on the module $(H\R)_X$.	
\end{proposition}
\proof The condition \eqref{xdefnr} is fulfilled by all elements of $(H\R)_X(k_+)$ since it involves sums on subsets of $k_+$. Conversely if $a\in \R^k$ fulfills \eqref{xdefnr} and $Z\subset k_+$ let $$Z_+:=\{j\in Z\mid a_j>0\}, \ \ Z_-:=\{j\in Z\mid a_j<0\}$$ 
One has $0\leq \sum_{Z_+} a_j \leq \lambda$, $0\geq \sum_{Z_-} a_j \geq -\lambda$
and thus $-\lambda\leq \sum_{Z} a_j\leq \lambda$.\newline
To prove the second statement, let $Y=k_+,Y'=k'_+$ be finite pointed sets and consider the map given by the product
$$
m:\Vert H\R\Vert_1(Y)\wedge (H\R)_X(Y')\to (H\R)(Y\wedge Y')
$$
It associates to $(\alpha_i)\in \Vert H\R\Vert_1(Y)$, $\sum \vert \alpha_i\vert \leq 1$ and $(a_j)\in (H\R)_X(Y')$ the doubly indexed $b:=(b_{i,j})$, $b_{i,j}=\alpha_i a_j$ and one needs to show that $b\in (H\R)_X(Y\wedge Y')$. Let 
$$
Y_+=\{i\in Y\mid \alpha_i>0\}, \ Y_-=\{i\in Y\mid \alpha_i<0\},\ Y'_+=\{j\in Y'\mid a_j>0\},\ Y'_-=\{j\in Y'\mid a_j<0\}
$$
By the rule of signs the pairs $(i,j)$ for which $b_{i,j}>0$ form the union $Y_+\times Y'_+\cup Y_-\times Y'_-$ so that one gets 
$$
\sum_{b_{i,j}>0}b_{i,j}=\sum_{Y_+\times Y'_+}\alpha_i a_j+\sum_{Y_-\times Y'_-}(-\alpha_i)(- a_j)=\sum_{Y_+}\alpha_i\sum_{Y'_+} a_j+\sum_{Y_-}(-\alpha_i)\sum_{ Y'_-}(- a_j)\leq \lambda
$$
using \eqref{xdefnr} for the sums over the $a_j$ together with the inequality  
$\sum_{Y_+}\alpha_i+\sum_{Y_-}(-\alpha_i)\leq 1$ (since $\sum \vert \alpha_i\vert \leq 1$). One treats in a similar way the sum over the negative $b_{i,j}$.\endproof 

In general, let $\sigma\in \Hom_\gop(k_+,1_+)$ with $\sigma(\ell)=1$ $\forall \ell\neq *$ and
$
\delta(j,k)\in \Hom_\gop(k_+,1_+)$, $\delta(j,k)(\ell):= 1$ if $\ell=j$, $\delta(j,k)(\ell):= *$ if $\ell\neq j$. \newline
Given an $\sss$-module $\cF$ and  elements $x,x_j\in \cF(1_+)$, $j=1,\ldots, k$, one writes 
\begin{equation}\label{defnsum}
x=\sum_j x_j\iff \exists z	\in \cF(k_+)~\text{s.t.}~  \cF(\sigma)(z)=x,~ \cF(\delta(j,k))(z)=x_j,\ \forall j.
\end{equation}
A tolerance relation $\cR$ on a set $X$ is a reflexive and symmetric relation on  $X$.  Equivalently, $\cR$ is a subset $\cR\subset X\times X$ which is symmetric and contains the diagonal. We shall denote by $\cT$ the category of tolerance relations $(X,\cR)$. Morphisms in  $\cT$  are defined by \[ 
\Hom_\cT((X,\cR),(X',\cR')):=\{	\phi: X \to X', \  \phi(\cR)\subset \cR'\}.
\]
We denote  $\cT_*$ the pointed category under the object $\{\ast\}$ endowed with the trivial relation. A tolerant $\sss$-module is  a pointed covariant functor $\gop \longrightarrow \cT_*$ (\cite{RR}). We recall below the definition of their dimension.
\begin{definition} \cite{RR}\label{generator}
 Let $(E,\cR)$ be a tolerant $\sss$-module. 
	A subset $F\subset E(1_+)$ generates $E(1_+)$ if the following two conditions hold
	\begin{enumerate}
	\item For $x,y\in F$, with $x\neq y ~\Longrightarrow~ (x,y)\notin \cR$	
	\item For every   $x\in E(1_+)$ there exists  $\alpha_j\in \{0,1\}$,  $j\in F$ and $y\in E(1_+)$ such that $y=\sum_F \alpha_j j\in E(1_+)$ in the sense of \eqref{defnsum},  and $(x,y)\in \cR$.
	\end{enumerate}
  The dimension $\dim_{\sss}(E,\cR)$ is defined as the minimal cardinality of a generating set $F$.
	\end{definition}

\section{Dimension of $H^0$ over $\sss$}

Let $m\in \N$, and $I_m=[-m,m]\cap \Z$. Next lemma follows from \eqref{defnsum} and Definition \ref{generator}.

\begin{lemma} \label{sgenerator} The dimension $\dim_\sss((\Z)_{I_m})$ is the smallest cardinality of a subset $G\subset I_m$ such that for any $j\in I_m$ there exists a subset $Z\subset G$ with $\sum_Z i=j$ and  $\sum_{Z'} i\in I_m$ for any $Z'\subset Z$.
\end{lemma}
The number of elements of $I_m$ is $2m+1$ and the number of subsets of $G$ is $2^{\# G}$, thus one has the basic inequalities 
\begin{equation}\label{basicinequ}
\# G\geq \log_2(2m+1)>\log_2(2m), \qquad \dim_\sss((\Z)_{I_m}\geq \lceil \log_2(m)\rceil +1.
\end{equation}
Here  $x\mapsto \lceil x\rceil$ denotes the ceiling function which associates to $x$ the smallest integer $> x$. For
$m=1$ one needs the two elements $\{-1,1\}$ to generate, while for $m=2$ one selects the three elements 
$\{-2,1,2\}$. For $m=3$ one takes the three elements 
$\{-3,1,2\}$ while for $m=4$ one takes the $4$ elements $\{-3,-1,1,3\}$.\vspace{.03in}

In general, one uses the following result.
\begin{lemma}\label{sequencebleft} Let $n\in \N$ and $I:=[-a,a]\subset \Z$, where $2^{n-1}\leq a<2^n$. \newline
$(i)$~If $n>4$ there exist $n$ distinct elements $\alpha_j\in (0,a)$ such that $\sum \alpha_j=a$ and that any element $z\in [0,a]$ can be written as a partial sum $z=\sum_Z\alpha_j$.\newline
$(ii)$~The minimal number of $\sss$-generators of $(H\Z)_I$ is $n+1$.  	
\end{lemma}
\proof $(i)$~We have $\sum_0^{n-1} 2^j=2^n-1\geq a$ and $\sigma:=\sum_0^{n-2} 2^j=2^{n-1}-1< a$. The idea is to adjoin to the set $T:=\{2^j\mid 0\leq j\leq n-2\}$, whose cardinality is $n-1$ and whose sum is $\sigma<a$, another element $a-\sigma$ so that the full sum is $a$. The first try is by taking $F=T\cup \{a-\sigma\}$. Assume first that $a-\sigma\notin T$. The partial sums obtained from $F$ are the union of the interval $[0,\sigma]$ with the interval $[a-\sigma,a]$ and these two intervals cover $[0,a]$, since $a-\sigma+\sigma=a$ while $a-\sigma\leq \sigma+1$. If $a-\sigma\in T$ one has for some $k\geq 0$ that $a=\sigma+2^k$. To avoid the repetition we adopt the following rules for $2^{n-1}\leq a<2^{n}$
\begin{enumerate}
\item If $a=2^{n-1}$ we let $F:=\{2^j\mid 0\leq j\leq n-3\}\cup \{2^{n-2}-2\}\cup \{3\}$
\item If  $a\neq 2^{n-1}$ and $a-\sigma\in T$, let $F:=\{2^j\mid 0\leq j\leq n-3\}\cup \{2^{n-2}-1\}\cup \{a-\sigma+1\}$
\item If  $a\neq 2^{n-1}$ and $a-\sigma\notin T$, let $F:=T\cup \{a-\sigma\}$
\end{enumerate}
Since by hypothesis $n>4$ one has $2^{n-2}-2>2^{n-3}$, so in case 1. one gets $\#F=n$ and the sum of elements of $F$ is $a=2^{n-1}$. The partial sums of elements of $\{2^j\mid 0\leq j\leq n-3\}$ cover the interval $J=[0,2^{n-2}-1]$. By adding $2^{n-2}-2$ to elements of $J$ one obtains the interval $J+2^{n-2}-2=[2^{n-2}-2,2^{n-1}-3]$ whose union with $J$ is $[0,2^{n-1}-3]$, then by imputing  the element $3\in F$ one sees that the partial sums cover $[0,a]$.\newline
In case 2. one obtains similarly  $\#F=n$ since $a-\sigma+1\notin T$  and the sum of elements of $F$ is $\sigma+a-\sigma=a$. The partial sums of elements of $\{2^j\mid 0\leq j\leq n-3\}$ cover the interval $J=[0,2^{n-2}-1]$ and using $2^{n-2}-1$ added to elements of $J$ one obtains the interval $J+2^{n-2}-1=[2^{n-2}-1,2^{n-1}-2]$ whose union with $J$ is $J'=[0,2^{n-1}-2]=[0,\sigma-1]$. Adding $a-\sigma+1$ to $J'$ one obtains the interval $J"=[a-\sigma+1,a]$. Since $a-\sigma\in T$ one has $a-\sigma\leq 2^{n-2}$, hence $a-\sigma+1\leq \sigma-1$, so that the lowest element of $J"$ belongs to $J'$ and $J'\cup J"=[0,a]$.\newline
In case 3. the partial sums of elements of $F$ cover $[0,a]$ as explained above.\newline
$(ii)$~Let $k$ be the minimal number of $\sss$-generators of $(H\Z)_I$. By \eqref{basicinequ} one has $k\geq n+1$. It remains to show that there exists a generating set of cardinality $n+1$. We assume first that $n>4$ and thus,  by $(i)$, let  $\alpha_j\in (0,a)$ be $n$ distinct elements fulfilling $(i)$. Let $F=\{-a\}\cup\{\alpha_j\}\subset [-a,a]$. By construction $\#F=n+1$. To show that $F$ is an $\sss$-generating set of $(H\Z)_I$ one needs to check the conditions of Lemma \ref{sgenerator}. By construction, the sum of positive elements of $F$ is $a$ and the sum of its negative elements is $-a$ thus any partial sum of elements of $F$ belongs to $I=[-a,a]$. Moreover the partial sums of positive elements of $F$ cover the interval $[0,a]$ by $(i)$, and using the element $-a$ one covers $I=[-a,a]$.\newline
For $n\leq 4$ one has $a\leq 15$ and one can list generating sets of cardinality $n+1$ as follows 
\begin{align*}
&\{-1,1\},\{-3,1,2\},\{-6,1,2,3\},\{-7,1,2,4\},\{-10,1,2,3,4\},\{-11,1,2,3,5\}\\ &\{-12,1,2,3,6\},\{-13,1,2,3,7\},\{-14,1,2,4,7\},\{-15,1,2,4,8\}
\end{align*}
These sets are of the same type as those constructed for $n>4$;  for the other values one has 
$$
\{-3,-1,1,3\},\{-4,-1,2,3\},\{-7,-1,1,2,5\},\{-8,-1,1,3,5\}.
$$
The value $a=2$ requires $3$ generators $\{-2,1,2\}$ and it is the only one for which the set $F$ of generators cannot be chosen in such a way that the sum of its positive elements is $a$ and the sum of its negative elements is $-a$. One nevertheless checks that all elements are obtained as an admissible sum. 
\endproof

\begin{theorem}\label{RRZ} Let $D$ be an Arakelov divisor  on $\spzb$. If $\deg(D)\geq 0$ one has
\begin{equation}\label{rrforz}
\dim_{\sss}H^0(D)=\bigg\lceil \deg_2 D\bigg\rceil	+1.
\end{equation}	
\end{theorem}
\proof One may assume that $D=\delta\{\infty\}$ where $\delta=\deg(D)>0$. One has $H^0(D)=(H\Z)_I$ where $I=[-e^\delta,e^\delta]$, using the classical relation between the degree of the divisor and the associated compact subset in adeles\footnote{Note that $e^{\deg D}=2^{\deg_2(D)}$}. Let $n\in \N$, $n\geq 1$, such that $2^{n-1}\leq e^\delta<2^n$. The integer part $a$ of $e^\delta$ fulfills $2^{n-1}\leq a<2^n$ and one has $H^0(D)=(H\Z)_{[-a,a]}$. Thus, by Lemma \ref{sequencebleft} one gets 
$\dim_{\sss}H^0(D)=n+1$. By definition  $\deg_2 D:=\deg D/\log 2$. The conditions $2^{n-1}\leq e^\delta<2^n$ mean that $n-1\leq \deg_2 D<n$ and show that the  least integer $>\deg_2 D$ is equal to $n$ which proves \eqref{rrforz}.\endproof 

\section{Dimension of $H^1$ over $\sss$}
We define the following sequence of integers:
\begin{equation}\label{bdefn}
 j(n):=\frac 13 (-2)^{n}-\frac 12
(-1)^{n}+\frac 16\qquad n\in\N.\end{equation}
 The first values of $j(n)$ are then:
$
0,1,-2,5,-10,21,-42,85,-170,341,-682,1365,-2730,\ldots
$
\begin{lemma}\label{sequenceb} Let  $G(n)=\{(-2)^j\mid 0\leq j<n\}$.  The map $\sigma$ from the set of subsets of $G(n)$ to $\Z$ defined by $\sigma(Z):=\sum_Z j$ is a bijection with the interval $\Delta(n):=[j(k),j(k)+2^n-1]$ 	where $k=k(n):=2E(n/2)+1$,  ($E(x)=$  integral part of $x$).	
\end{lemma}
\proof The map $\sigma$ is injective and covers an interval $[a,b]$. The lower bound $a$ is the sum of powers $a=\sum_{0\leq \ell <\frac{n-1}{2}}(-2)^{2\ell+1}$ and the upper bound is the sum of powers $b=\sum_{0\leq \ell <\frac{n}{2}}(-2)^{2\ell}$. We list the first intervals  as follows
$$
\Delta(1)=[0,1], \ \ \Delta(2)=[-2,1], \ \ \Delta(3)=[-2,5], \ \ \Delta(4)=[-10,5], \ \ \Delta(4)=[-10,21],\ldots 
$$
\endproof
 We refer to \cite{RR}, Appendix A, B, for the interpretation of $H^1(D)$ in terms of the tolerant $\sss$-module  $(U(1),d)_\lambda$, $\lambda=e^{\deg D}$. At level $1$ the tolerance relation on the abelian group $\R/\Z$ is given by the condition $d(x,y)\leq \lambda$.  
\begin{proposition}\label{caseh1} Let  $U(1)$ be the abelian group $\R/\Z$ endowed with the canonical metric $d$ of length $1$.
 Let $\lambda\in\R_{>0}$, $U(1)_\lambda$ the tolerant $\sss$-module $(U(1),d)_\lambda$. 	Then  
 \begin{equation}\label{dimh1}
\dim_{\sss}U(1)_\lambda=\begin{cases}m  &\text{if\  $2^{-m-1}\leq\lambda< 2^{-m}$,}\\0 &\text{if $\lambda\geq \frac 12$.}	
\end{cases}
\end{equation}
 \end{proposition}
 \proof For $\lambda\geq \frac 12$, any element of $U(1)_\lambda=(\R/\Z,d)_\lambda$ is at distance $\leq \lambda$ from $0$, thus one can take $F=\emptyset$ as generating set since, by convention, $\sum_\emptyset=0$. Thus $\dim_{\sss}U(1)_\lambda=0$. Next, we assume $\lambda<\frac 12$. Let $F\subset U(1)$ be a generating set and let $k=\# F$. One easily sees that there are at most $2^k$ elements of the form $\sum_F \alpha_j j$, $\alpha_j\in \{0,1\}$. The subsets  $\{x\in U(1) \mid d(x,\sum_F\alpha_j j)\leq \lambda\}$  cover $U(1)$, and since each of them has measure $2\lambda$ one gets the inequality $2 \lambda\cdot  2^k\geq 1$. Thus  $k\geq \frac{-\log \lambda- \log 2}{\log 2}$. When $\frac{-\log \lambda- \log 2}{\log 2}=m$ is an integer, one has $\lambda= 2^{-m-1}$. Let $F(m)=\{(-2)^{-j}\mid 1\leq j\leq m\}$. The minimal distance between two elements of $F(m)$ is the distance between $2^{-m+1}$ and $-2^{-m}$ which is $3 \cdot 2^{-m}=6 \lambda$.
Let us show that $F(m)$ is a generating set.   By Lemma \ref{sequenceb} any integer $q$ in the interval $\Delta(m)$ can be written as  $q=\sum_{i=0}^{m-1} \alpha_i (-2)^i$, with $\alpha_i\in\{0,1\}$. One then gets 
\[
q \cdot (-2)^{-m}=\sum_{i=0}^{m-1} \alpha_i (-2)^{i-m}=\sum_{j=1}^{m}\alpha_{m-j}(-2)^{-j}.
\]
 Let $y\in \R/\Z$, lift $y$ to an element $x$ of the interval $(-2)^{-m}[j(k(m)),j(k(m))+2^m)$ which is connected of length $1$ and is a fundamental domain for the action of $\Z$ by translation.  
 Then  there exists an integer $q\in \Delta(m)$ such that $ \vert (-2)^m  x-q\vert \leq \frac 12$. 
 Hence $d(x,q \cdot (-2)^{-m})\leq 2^{-m-1}=\lambda$. This proves that $F(m)$ is a generating set (see  Definition \ref{generator}) and one derives $\dim_{\sss}U(1)_\lambda=m$. Assume now that $\frac{-\log \lambda- \log 2}{\log 2}\in(m,m+1)$, where $m$ is an integer, \ie that $\lambda\in (2^{-m-2},2^{-m-1})$. For any generating set $F$ of cardinality $k$ one has $k\geq \frac{-\log \lambda- \log 2}{\log 2}>m$ so that $k\geq m+1$. The subset $F(m+1)=\{(-2)^{-j}\mid 1\leq j\leq m+1\}$ fulfills the first condition of Definition \ref{generator} since the minimal distance between two elements of  $F(m+1)$ is $3 \cdot 2^{-m-1}$ which is larger than $\lambda< 2^{-m-1}$. As shown above, the subset $F(m+1)$ is generating for $\lambda=2^{-m-2}$ and a fortiori for $\lambda>2^{-m-2}$. Thus one obtains $\dim_{\sss}U(1)_\lambda=m+1$ and  \eqref{dimh1} is proven.\endproof 
 \section{Riemann-Roch formula}
\begin{figure}[H]	\begin{center}
\includegraphics[scale=0.45]{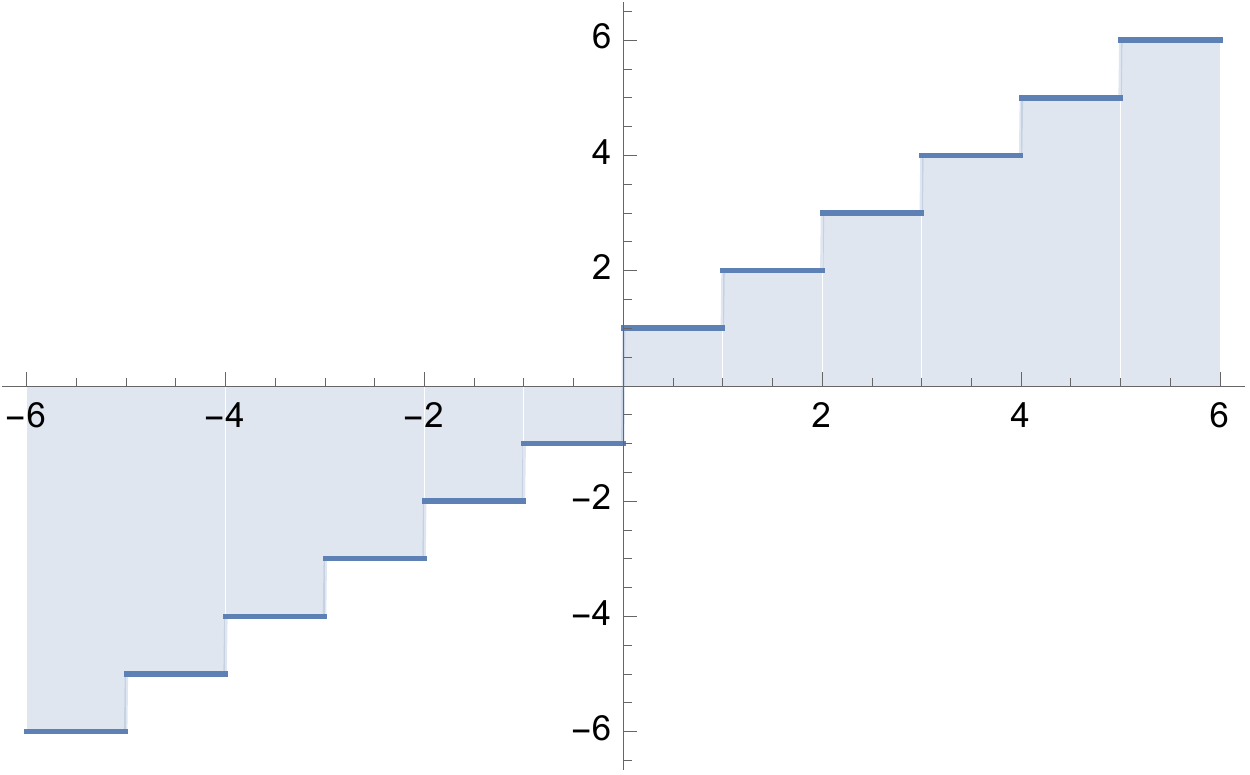}
\end{center}
\caption{Graph of $\dim_{\sss}H^0(D)-\dim_{\sss}H^1(D)-1$ as a function of $\deg_2 D$ \label{rr1}}
\end{figure}

 We can now formulate the main result of our paper
 \begin{theorem}\label{rrspzb} 
Let $D$ be an Arakelov divisor  on $\spzb$. Then 
\begin{equation}\label{rrforq1}
\dim_{\sss}H^0(D)-\dim_{\sss}H^1(D)=\bigg\lceil \deg_2 D\bigg\rceil'	+1
\end{equation}
where $\lceil x \rceil'$ is the right continuous function which agrees with ceiling$(x)$ for $x>0$ non-integer and with -ceiling$(-x)$ for $x<0$ non-integer (see Figure \ref{rr1}).
\end{theorem}
\proof For $\deg_2 D\geq 0$ one has $\lambda=e^{\deg D}\geq 1$
and hence by \eqref{dimh1} one gets $\dim_{\sss}H^1(D)=0$, so 
\eqref{rrforq1} follows from Theorem \ref{RRZ}. For $\deg_2 D< 0$ one has $\dim_{\sss}H^0(D)=0$ since the empty set is a generating set. For $\deg_2 D\in [-m-1,-m)$ where $m\in \N$ one has, by \eqref{dimh1},  $\dim_{\sss}H^1(D)=m$. Thus the left hand side of \eqref{rrforq1} is $-m$ while the right hand side is 
$$
\bigg\lceil \deg_2 D\bigg\rceil'	+1=-m
$$
by definition of the function $\lceil x \rceil'$ as the right continuous function which agrees with ceiling$(x)$ for $x>0$ non-integer and with -ceiling$(-x)$ for $x<0$ non-integer.\endproof

\begin{Backmatter}

\paragraph{Acknowledgments}
The second author is partially supported by the Simons Foundation collaboration grant n. 691493.

\end{Backmatter}

\end{document}